\documentclass{gtart_h}

\def\ifplaintex{\expandafter\ifx\csname documentclass\endcsname\relax}

\ifplaintex 
\else
\fi

\expandafter\ifx\csname beginpicture\endcsname\relax
\expandafter\ifx\csname documentclass\endcsname\relax
\input pictex \else
\input prepictex \input pictex \input postpictex \fi\fi

\def\gt{{\mathsurround=0pt\it $\cal G\mskip-2mu$eometry \&\ 
$\cal T\!\!$opology}}        

\def\gtp{{\mathsurround=0pt\it $\cal G\mskip-2mu$eometry \&\ 
$\cal T\!\!$opology $\cal P\!$ublications}}  

\def\lognumber#1{\def\thelognumber{#1}}
\def\volumenumber#1{\def\thevolumenumber{#1}}
\def\papernumber#1{\def\thepapernumber{#1}}
\def\volumeyear#1{\def\thevolumeyear{#1}}

\def\pagenumbers#1#2{\def\startpage{#1}\def\finishpage{#2}}
\def\published#1{\def\publishdate{#1}}
\def\proposed#1{\def\theproposer{#1}}
\def\seconded#1{\def\theseconders{#1}}
\def\received#1{\def\receiveddate{#1}}
\def\revised#1{\def\reviseddate{#1}}
\def\accepted#1{\def\accepteddate{#1}}

\long\def\asciiabstract#1{\long\def\theasciiabstract{#1}}

\let\\\par\let\thelognumber\relax
\let\thevolumenumber\relax\let\thepapernumber\relax
\let\thevolumeyear\relax\let\thesamplenumber\relax\let\startpage\relax
\let\finishpage\relax\let\publishdate\relax\let\receiveddate\relax
\let\reviseddate\relax\let\accepteddate\relax\let\theasciititle\relax
\let\theasciiauthors\relax
\let\theasciiabstract\relax
\let\theasciiemail\relax\let\theshortauthors\relax\let\theshorttitle\relax

\long\def\maketitlep{   

\count0=\startpage

\gt\hfill      
\beginpicture
\setcoordinatesystem units <0.33truein, 0.33truein> point at 2.2 0.9
\setplotsymbol ({$\cal G$})
\plotsymbolspacing=9truept
\circulararc 315 degrees from 0 1 center at 0 0
\setplotsymbol ({$\cal T$})
\circulararc 315 degrees from 1 -1 center at 1 0
\endpicture
\break
{\small\ifx\thesamplenumber\relax 
Volume \else Sample
\fi\thevolumenumber\ (\thevolumeyear)
\startpage--\finishpage\nl
Published: \publishdate}
\vglue 0.5truein plus 0.4fil minus 0.1truein

{\parskip=0pt\leftskip 0pt plus 1fil\def\\{\par\smallskip}{\ifplaintex\large
\else\Large\fi\bf\thetitle}\par\medskip}   

\vglue 0pt plus 0.1fil 

{\parskip=0pt\leftskip 0pt plus 1fil\def\\{\par}{\sc\theauthors}
\par\medskip}

\vglue 0pt plus 0.1fil 

{\small\parskip=0pt\let\newline\\
{\leftskip 0pt plus 1fil\def\\{\par}{\sl\theaddress}\par}
\expandafter\ifx\theemail\relax    
\relax\else\vglue 5pt plus 0.02fil minus 2pt\def\\{\stdspace{\rm 
and}\stdspace} 
\cl{Email:\stdspace\tt\theemail}\fi
\ifx\theurl\relax                  
\relax\else\vglue 5pt plus 0.02fil minus 2pt\def\\{\stdspace{\rm 
and}\stdspace}
\cl{URL:\stdspace\tt\theurl}\fi\par}

\vglue 7pt plus 0.3fil minus 3pt

{\bf Abstract}
\vglue 5pt plus 0.1fil minus 2pt

\theabstract

\vglue 7pt plus 0.3fil minus 3pt

{\bf AMS Classification numbers}\quad Primary:\quad \theprimaryclass

Secondary:\quad \thesecondaryclass

\vglue 5pt plus 0.3fil minus 2pt

{\bf Keywords:}\quad \thekeywords

\vglue 10pt plus 0.5fil minus 5pt

{\small  Proposed: \theproposer\hfill Received: \receiveddate\nl
Seconded: \theseconders\hfill 
\ifx\reviseddate\relax                         
Accepted: \accepteddate                        
\else
Revised: \reviseddate                          
\fi}
\eject
}       

\let\maketitlepage\maketitlep
\let\maketitle\maketitlepage

\font\phead=cmsl9 scaled 950
\font=cmsl9 scaled 1050
\font\pnum=cmbx10 scaled 913
\font\lnum=cmbx10 
\font\pfoot=cmsl9 scaled 950
\font=cmsl9 scaled 1050
\ifplaintex
\headline{\vbox to 0pt{\vskip -4.5mm\line{\small\phead\ifnum
\count0=\startpage ISSN 1364-0380 (on line)
1465-3060 (printed) \hfill {\pnum\folio}\else\ifodd\count0\def\\{ }%
\ifx\theshorttitle\relax\thetitle\else\theshorttitle\fi\hfill{\pnum\folio}
\else\def\\{ and }{\pnum\folio}\hfill\ifx\theshortauthors\relax\theauthors
\else\theshortauthors\fi\fi\fi}\vss}}
\footline{\vbox to 0pt{\vglue 0mm\line{\small\pfoot\ifnum\count0=\startpage
\copyright\ \gtp\hfill\else
\gt, Volume \thevolumenumber\ (\thevolumeyear)\hfill\fi}\vss
}}
\else
\makeatletter
\def\@oddhead{{\small\lhead\ifnum\count0=\startpage ISSN 1364-0380 (on line)
1465-3060 (printed) \hfill {\lnum\number\count0}\else\ifodd\count0
\def\\{ }\ifx\theshorttitle\relax \thetitle \else\theshorttitle\fi\hfill
{\lnum\number\count0}\else\def\\{ and }{\lnum\number\count0}
\hfill\ifx\theshortauthors\relax 
\theauthors\else\theshortauthors\fi\fi\fi}}\def\@evenhead{\@oddhead}
\def\@oddfoot{\small\lfoot\ifnum\count0=\startpage\copyright\ \gtp\hfill\else
\gt, Volume \thevolumenumber\ (\thevolumeyear)\hfill\fi}
\def\@evenfoot{\@oddfoot}
\makeatother
\fi

\newwrite\gtoutfile
\long\gdef\makeheadfile{  
{\def\\{, }\def\s{ }
\immediate\openout\gtoutfile head.xxx
\immediate\write\gtoutfile{Proxy-for: \ifx\theasciiauthors\relax
\theauthors\else\theasciiauthors\fi\s<\ifx\theasciiemail\relax\theemail\else\theasciiemail\fi>}
\immediate\write\gtoutfile{\noexpand\\}
\immediate\write\gtoutfile{Authors: \ifx\theasciiauthors\relax
\theauthors\else\theasciiauthors\fi}
{\def\\{ }\immediate\write\gtoutfile{Title: \ifx\theasciititle\relax
\thetitle\else\theasciititle\fi}}
\immediate\write\gtoutfile{Subj-class: GT or SG or MG etc}
\immediate\write\gtoutfile{MSC-class: \theprimaryclass\ifx\thesecondaryclass\relax\else, \thesecondaryclass\fi}
\immediate\write\gtoutfile{Journal-ref: Geom. Topol. \thevolumenumber
(\thevolumeyear) \startpage-\finishpage}
\immediate\write\gtoutfile{Comments: Published by Geometry and Topology at}
\immediate\write\gtoutfile{\s\s http://www.maths.warwick.ac.uk/gt/GTVol\thevolumenumber/paper\thepapernumber.abs.html}
\immediate\write\gtoutfile{\noexpand\\}
\immediate\write\gtoutfile{}
\ifx\theasciiabstract\relax
\immediate\write\gtoutfile{\theabstract}\else
\immediate\write\gtoutfile{\theasciiabstract}\fi
\immediate\write\gtoutfile{}
\immediate\write\gtoutfile{\noexpand\\}
\immediate\write\gtoutfile{}
\immediate\closeout\gtoutfile}}  

\def\maketitlepage{\maketitlep\makeheadfile}
\let\maketitle\maketitlepage

\lognumber{384}
\volumenumber{8}\papernumber{6}\volumeyear{2004}
\pagenumbers{277}{293}
\received{25 November 2003}
\revised{13 January 2004}
\published{14 February 2004%
\vglue 3mm\sl To Ada%
}
\accepted{2 January 2004}
\proposed{Leonid Polterovich}
\seconded{Peter Ozsv\'ath, Dieter Kotschick}

\usepackage{amsmath, amssymb}

\def\e{\varepsilon}

\def\g{\gamma}

\newcommand{\Flux}{\mathrm{Flux}}
\newcommand{\Hol}{\mathrm{Hol}}

\newcommand\hook{\mathbin{\hbox{\vrule height .5pt width 3.5pt depth 0pt
\vrule height 6pt width .5pt depth 0pt}}}

\newcommand{\wt}{\widetilde}
\newcommand{\wh}{\widehat}

\newcommand{\cD}{\mathcal {D}}

\newcommand{\Ker}{{\mathrm {Ker}}}

  \newcommand{\genus}{{\mathrm {genus}}}

\newcommand{\Id}{{\mathrm {Id}}}

\newcommand{\D}{{\rm D}}

\newcommand{\Z}{{\mathbb Z}}

\newcommand{\R}{\mathbb {R}}

\newcommand{\C}{\mathbb {C}\,}

\newcommand{\const}{\mathrm {const}}

\newcommand{\no}{\noindent}

\newtheorem{theorem}{Theorem}[section]
\newtheorem{corollary}[theorem]{Corollary}
\newtheorem{lemma}[theorem]{Lemma}
\newtheorem{proposition}[theorem]{Proposition}

\newtheorem{remark}[theorem]{Remark}

\def\R{\mathbb{R}}
\def\Z{\mathbb{Z}}

\def\cD{{\cal D}}

\begin{document}

\title{A few remarks about symplectic filling}
\author{Yakov Eliashberg}
\address{Department of Mathematics,
Stanford University\\Stanford CA 94305-2125,
USA}
\email{eliash@gauss.stanford.edu}

\begin{abstract}
We show that any compact symplectic manifold $(W,\omega)$ with
boundary embeds as a domain into a closed symplectic manifold,
provided that there exists a contact plane $\xi$ on $\partial W$ which
is weakly compatible with $\omega$, i.e.\ the restriction
$\omega|_\xi$ does not vanish and the contact orientation of $\partial
W$ and its orientation as the boundary of the symplectic manifold $W$
coincide.  This result provides a useful tool for new applications by
Ozsv\'ath--Szab\'o of Seiberg--Witten Floer homology theories in
three-dimensional topology and has helped complete the
Kronheimer--Mrowka proof of Property P for knots.
\end{abstract}

\asciiabstract{We show that any compact symplectic manifold (W,\omega)
with boundary embeds as a domain into a closed symplectic manifold,
provided that there exists a contact plane \xi on dW which is weakly
compatible with omega, i.e. the restriction \omega|\xi does not vanish
and the contact orientation of dW and its orientation as the boundary
of the symplectic manifold W coincide.  This result provides a useful
tool for new applications by Ozsvath-Szabo of Seiberg-Witten Floer
homology theories in three-dimensional topology and has helped
complete the Kronheimer-Mrowka proof of Property P for knots.}

\keywords{Contact manifold, symplectic filling, symplectic Lefschetz
fibration, open book decomposition}

\primaryclass{53C15}
\secondaryclass{57M50}

\maketitle

\section{Introduction}
All manifolds which we consider in this article are assumed
oriented. A contact manifold $V$ of dimension three carries a
canonical orientation. In this case we will denote by $-V$ the contact
manifold with the opposite orientation. Contact plane fields are
assumed co-oriented, and therefore oriented.  Symplectic manifolds are
canonically oriented, and so are their boundaries.

\no We prove in this article the following theorem:
\begin{theorem}\label{thm:cobord}
Let $(V,\xi)$ be a contact manifold and $\omega$ a closed $2$--form on
$V$ such that $\omega|_{\xi}>0$. Suppose that we are given an open
book decomposition of $V$ with a binding $B$.  Let $V'$ be obtained
from $V$ by a Morse surgery along $B$ with a canonical $0$--framing, so
that $V'$ is fibered over $S^1$. Let $W$ be the corresponding
cobordism, $\partial W=(-V)\cup V'$. Then $W$ admits a symplectic form
$\Omega$ such that $\Omega|_V=\omega$ and $\Omega$ is positive on
fibers of the fibration $V'\to S^1$.
\end{theorem}

\begin{remark} {\rm Note that the binding $B$ has a canonical decomposition of its tubular neighborhood
given by the pages of the book. The $0$--surgery along $B$ is the Morse surgery  associated
with this decomposition. If the binding is disconnected then we assume that the surgery is performed simultaneously
along all the components of $B$.}
\end{remark}

\no We will deduce the following result from Theorem \ref{thm:cobord}:
\begin{theorem}\label{thm:concave-filling}
Let  $(V,\xi)$ and $\omega$ be as  in Theorem \ref{thm:cobord}.
Then there exists a symplectic manifold $(W',\Omega')$ such that $\partial W'=-V$ and  $\Omega'|_V=\omega$.
Moreover, one can arrange that
$H_1(W')=0$,\footnote{This observation is due to Kronheimer and Mrowka, see \cite{Kronheimer-Mrowka}.} and
that $(W',\Omega')$ contains the symplectic cobordism $(W, \Omega)$ constructed
in Theorem \ref{thm:cobord} as a subdomain adjacent to the boundary.
In particular,  any symplectic manifold  which  weakly fills
(see Section \ref{sec:flavors} below) the contact manifold
$(V,\xi)$ can be symplectically embedded  as a subdomain into a closed symplectic manifold.
\end{theorem}

\begin{corollary}\label{cor:semi}
Any  weakly (resp.\ strongly) semi-fillable (see \cite{Eliash-Thurston})
contact manifold is weakly (resp.\ strongly) fillable.
 \end{corollary}
 \begin{remark}{\rm
  Theorem \ref{thm:cobord} serves as a missing  ingredient in proving that the Ozsv\'ath--Szab\'o contact invariant
 $c(\xi)$ does not vanish for weakly symplectically fillable (and hence for non-existing anymore semi-fillable)
 contact structures. This and other applications of the results of this article in the Heegaard Floer homology theory
are discussed
  in the  paper of Peter Ozsv\'ath and Zolt\'an Szab\'o, see \cite{Ozsvath-Szabo}.
 The observation made in this paper
  also helped to streamline the  program of Peter Kronheimer and  Tomasz Mrowka for proving the Property P for knots,
see their  paper \cite{Kronheimer-Mrowka}.}
  \end{remark}

 \noindent{\bf Acknowledgements}\qua This article is my answer to the
 question I was asked by Olga Plamenevskaya and David Gay during the
 workshop on Floer homology for $3$--manifolds in Banff International
 Research Station.  I want to thank them, as well as Selman Akbulut,
 Michael Freedman, Robion Kirby, Peter Kronheimer, Robert Lipshitz,
 Tomasz Mrowka, Peter Ozsv\'ath, Michael Sullivan and Zolt\'an Szab\'o
 for stimulating discussions, and Augustin Banyaga, Steven Kerckhoff,
 Leonid Polterovich and Andr\'as Stipsicz for providing me with the
 necessary information. I am grateful to John Etnyre, Tomasz Mrowka,
 Andr\'as Stipsicz and Leonid Polterovich for their critical remarks
 concerning a preliminary version of this paper. I also want to thank
 Peter Kronheimer and Tom Mrowka with sharing with me their
 alternative ``flux-fixing argument" (see Lemma \ref{lm:thimble}
 below).
 
This research was partially supported by NSF grants DMS-0204603 and
DMS-0244663.

\section{Proof of Theorem \ref{thm:cobord}}
 We begin with the following lemma which is a slight reformulation of Proposition 3.1 in \cite{Eliash-boundary}.
 A similar statement is contained also in \cite{Kron-Mrowka-monopole}.
 \begin{lemma}\label{lm:bound-behavior}
  Let $(V,\xi)$ and $\omega$ be as in Theorem \ref{thm:cobord}. Then given any contact form
 $\alpha $ for $\xi$ and any $C>0$ one can find
  a symplectic form $\Omega$ on $V\times[0,1]$ such that
  \begin{itemize}
\item[\rm a)] $\Omega|_{V\times 0}=\omega$;
\item[\rm b)] $\Omega_{V\times [1-\varepsilon,1]}=\omega +Cd(t\alpha)$, where $t\in[1-\varepsilon,1]$ and $0<\varepsilon<1$;
\item[\rm c)] $\Omega$ induces the negative orientation on $V\times 0$ and positive on $V\times 1$.
\end{itemize}
\end{lemma}
\begin{proof}
By assumption $$\omega|_\xi= fd\alpha|_\xi=d(f\alpha)|_\xi$$ for a positive function $V\to\R$.
Set $\tilde\alpha=f\alpha$. Then $\omega=d\tilde\alpha+\tilde\alpha\wedge\beta$.
Take a smooth  function $h\co V\times[0,1]\to\R$ such that $$h|_{V\times 0}=0,\;h|_{V\times [1-\e,1]}=\frac{Ct}{f},\;\;
\frac{dh}{dt}>0,$$ where $t$ is the coordinate corresponding to the projection $V\times[0,1]\to[0,1]$.
Consider the form
$\Omega =\omega+d(h\tilde\alpha)$. Here we keep the notation $\omega$ and $\tilde\alpha$ for the pull-backs of
$\omega$ and $\tilde\alpha$ to $V\times[0,1]$.
Then we have
\begin{equation*}
\Omega=d\tilde\alpha+\tilde\alpha\wedge\beta+d_Vh\wedge\tilde\alpha+\frac{dh}{dt}dt\wedge \tilde\alpha,
\end{equation*}
where $d_Vh$ denotes the differential of $h$ along $V$.
Then
\begin{equation*}
\Omega\wedge\Omega=2\frac{dh}{dt} dt\wedge \tilde\alpha\wedge d\tilde\alpha > 0\,.
\end{equation*}
Hence $\Omega$ is symplectic and it clearly satisfies the conditions a)--c).
\end{proof}

\no Let us recall that  a contact form $\lambda$
 on   $V$ is called {\it  compatible} with the given open book decomposition
 (see \cite{Giroux}) if
\begin{itemize}
\item[a)] there exists a neighborhood $U$ of the binding  $B$, and the  coordinates
($r,\varphi,u)\in[0,R]\times\R/{2\pi\Z}\times {\R/2\pi \Z}$ such that $$U=\{r\leq R\}\;\;\hbox{  and}\;\;
 \lambda|_U=h(r)(du+r^2d\varphi)\,,$$
where the positive  $C^\infty$--function $h$ satisfies the conditions $$h(r)-h(0)=-r^2\;
\hbox{ near }\;r=0\;\hbox{ and}\; h'(r)<0\;\hbox{for all}\; r>0\,;$$
\item[b)] the parts of pages of the book in $U$ are given by equations $\varphi=\const$;
\item[c)] $d\alpha$ does not vanish on the  pages of the book (with the binding deleted).
\end{itemize}

\begin{remark} {\rm
An admissible contact form $\alpha$ defines an orientation of pages and hence
an orientation of the binding $B$ as the boundary of a page. On the other hand, the form $\alpha$ defines
a co-orientation of the contact plane field, and hence an orientation of $B$ as a transversal curve.
These two orientations of $B$ coincide.
}
\end{remark}
\begin{remark}{\rm By varying admissible forms for a given contact plane field
one can arrange any function $h$ with the properties described in a).
Indeed, suppose we are  given another function $\wt h$ which satisfies  a). We can assume
without loss of generality that $\alpha$ has the presentation a) on a bigger domain $U'=\{r\leq R'\}$
for $R'>R$.
Let us choose $c>0$ such that $\wt h(R)>ch(R)$ and extend $\wt h$ to $[0,R']$ in such a way that
 $\wt h'(r)<0$ and $\wt h(r)=ch(r)$ near $R'$. Then the form $\wt h\alpha$ on $U'$ extended to the rest of the
manifold $V$ as $c\alpha$ is admissible for the given open book decomposition.
}
\end{remark}

\no Now we are ready to prove Theorem \ref{thm:cobord}.
\begin{proof}[Proof of Theorem \ref{thm:cobord}]
Take a  constant $a>0$ and consider a  smooth on $[0,1)$ function $g\co [0,1]\to \R$ such that
$g|_{[0,1/2]}=a$,  $g(t)=\sqrt{1-t^2}$ for $t$ near $1$ and  $g'<0$ on $(1/2,1)$.

\no In the standard symplectic $\R^4$ which we  identify with $\C^2$  with coordinates
$(z_1=r_1e^{i\varphi_1}=x_1+iy_1,z_2=r_2e^{i\varphi_2}=x_2+iy_2)$ let us consider a domain
\begin{equation*}
\wt P=\{r_1\leq g(r_2),\,r_2\in[0,1]\}\,.
\end{equation*}
The domain $\wt P$ is contained in the polydisc
 $P=\{r_1\leq a,r_2\leq 1\}$  and can be viewed as obtained by smoothing the corners of $P$.

\noindent Let us denote by $\Gamma$ the part of the boundary of $\wt P$ given by
\begin{equation*}
\Gamma=\{\{r_1=g(r_2),\,r_1\in[1/2,1]\}\,.
\end{equation*}
Note that $\Gamma$  is $C^\infty$--tangent to $\partial P$ near its boundary.
 The primitive
 $$\gamma=\frac12\left(r_1^2d\varphi_1+r_2^2d\varphi_2\right)\,$$
of  the standard symplectic form
$$\omega_0 =dx_1\wedge dy_1 +dx_2\wedge dy_2=r_1dr_1\wedge d\varphi_1+r_2dr_2\wedge d\varphi_2$$
restricts to $\Gamma$ as a contact form
\begin{equation*}
\gamma|_{\Gamma}=\frac{r_2^2}{2}\left(\frac{g^2(r_2)}{r_2^2} d\varphi_1+d\varphi_2\right)\,.
\end{equation*}
Consider the product $G=S^2\times\D^2$ with the split symplectic structure
$\omega_0=\sigma_1\oplus\sigma_2$, where the total area of the  form $\sigma_1$ on $S^2$ is equal to $2\pi $ and the
total area of the form $\sigma_2$ on the disc $D^2$ is equal to $\pi a^2$. Note that if $S^2_+$ and $S^2_-$ denote
the upper and lower hemispheres of $S^2$ of equal area, then there exists a symplectomorphism
$$\Phi\co P\to S^2_+\times D^2\subset S^2\times D^2=G\,.$$
Let  $H$ be the closure of $G\setminus \Phi(\wt P)$ and $\wt \Gamma$
denote the image $\Phi(\Gamma)\subset\partial H$. Note that $$\Delta=\overline{\partial H\setminus\Gamma}
=\wt S^2_-\times \partial D^2=\bigcup\limits_{x\in \partial D^2}\wt S^2_-\times x=
\bigcup\limits_{x\in \partial D^2}D_x  \,,$$
where $\wt S^2_-$, $S^2_-\subset \wt S^2_-\subset S^2$, is a disc of area $\frac{9\pi}{4}$.
Thus  $H$ is a $2$--handle whose boundary
$\partial H$ consists of $\Gamma$ and $\Delta$ which meets along an infinitely sharp corner $\Gamma\cap\Delta$.
The part $\Delta$ is fibered by discs $D_x,x\in\partial D^2$, which are symplectic with respect to the form $\omega_0$.

\noindent Consider now a contact form $\lambda$ on $V$ compatible with the given open book decomposition. In particular,
on a neighborhood $$U=[0,R]\times\R/{2\pi\Z}\times {\R/2\pi \Z}\subset V$$ of the binding $B$
we have
$\lambda|_U=h(r)(du+r^2d\varphi)$,
where the positive  $C^\infty$--function $h$ satisfies the conditions
$$h(r)-h(0)=-r^2\;\hbox{ near}\;\ r=0\;
 \hbox{and} \; h'(r)<0\;\hbox{ for all}\; r>0\,.$$
\noindent Let us choose $a=\frac R2$, and consider a diffeomorphism $F\co \Gamma\to U$ given by the formula
\begin{equation*}
r=\frac{g(r_2)}{r_2},\;\varphi=\varphi_1,\; u=\varphi_2\,.
\end{equation*}
The function $$r_2\mapsto\frac{g(r_2)}{r_2}$$ maps $[1,\frac12]$ onto $[0,2a]=[0,R]$.
Let $\psi\co [0,R]\to[1,\frac12]$ be the inverse function.
Then
\begin{equation*}
F_*\gamma=\frac{\psi^2(r)}2\left( r^2d\varphi+du\right)\,.
\end{equation*}
Hence $F$ is a contactomorphism
$$(\Gamma,\{\gamma=0\})\to (U,\{\lambda=0\}=\xi)\,.$$ Moreover,
the form $F_*\gamma$, extended  to $V$ as   $\lambda$ on $V\setminus U$, defines on $V$
a smooth contact form
compatible with the given open book decomposition.

\no
Now we use Lemma \ref{lm:bound-behavior} to define on
 the collar $V\times[0,1]$ a  symplectic form $\Omega$ which  satisfies the conditions \ref{lm:bound-behavior}a)--c),
 where the constant $C$ will be chosen later. In particular, near $V\times 1$  we have
 $\Omega=Cdt\lambda+\omega$.
 Viewing $\Gamma$ as a part of the boundary of the handle $H$, we can extend $F$  to a symplectomorphism, still
 denoted by $F$,
of a neighborhood of $\Gamma\subset H$ endowed with the standard symplectic structure $C\omega_0$ to a
neighborhood of $U\subset V= V\times 1$ in $V\times[0,1]$ endowed with the symplectic structure $Cd(t\lambda)$.
 Note that the closed form $F^*\omega$ is exact:
$$F^*\omega=d\theta,$$ and hence it extends to $H$ as $\wt \omega=d(\sigma\theta)$ where $\sigma$ is a cut-off function
equal to $0$ outside a neighborhood of $\Gamma$ in $H$. If $C$ is chosen sufficiently large
then the form $\Omega_0=C\omega_0+\wt\omega$ is symplectic, and its restrictions to the discs
$D_x\subset\Delta$, $x\in\partial D^2$, are symplectic as well. Hence the map $F$
can be used for  attaching
the symplectic handle $(H,C\omega_0+\wt\omega)$ to $V\times[0,1]$ along $U$.
The  resulted symplectic manifold $$W=V\times[0,1]\mathop{\cup}\limits_{U=F(\Gamma)}H$$ is the required
symplectic cobordism. Indeed we have
$$\partial (W,\Omega_0)=(-V,\omega)\cup(V',\omega'),$$ where
the component $V'$ of its boundary
is fibered over $S^1$ by closed  surfaces
 formed by parts of pages of the book inside $V\setminus U$ and discs
$D_x$. These surfaces are symplectic with respect to the form
$\omega'={\Omega_0}|_{V'}$.
\end{proof}

\section{Filling of symplectic fibrations over circle}\label{sec:fibrations}
A pair $(V,\omega)$, where $V$ is an oriented $3$--manifold  fibered over $S^1=\R/\Z$, and $\omega$
is  a closed $2$--form  which  is positive on the
fibers of the fibration, will be referred  to as a {\it symplectic fibration over $S^1$}. The projection
$V\to S^1$ will be denoted by $\pi$. We will assume that all symplectic fibrations we consider are normalized by the condition
that the integral of $\omega$ over a fiber is equal to $1$. The form $\omega$ induces a $1$--dimensional {\it characteristic }
foliation $\mathcal {F}_\omega$ on $V$ generated by the kernel of $\omega$. This foliation is
transversal
to the fibers of the fibration.
 The orientation of $V$ together with the symplectic orientation
of the fibers defines an orientation of $\mathcal {F}_\omega$.
Fixing a fiber $F_0$ over $0\in S^1=\R/\Z$  we can define  the {\it holonomy
diffeomorphism} $\Hol_{V,\omega}\co F_0\to F_0$. This is an area preserving diffeomorphism which defines $(V,\omega)$
uniquely up to a fiber preserving diffeomorphism fixed on $F_0$. Note that $\Hol_{-V,\omega}=\Hol^{-1}_{V,\omega}$.
Two symplectic fibrations are equivalent via an equivalence fixed on $F_0$ if and only if their
 holonomy diffeomorphisms coincide. If for
symplectic fibrations $(V,\omega_0)$ and $(V,\omega_1)$
  the holonomy diffeomorphisms $\Hol_{V,\omega_0}$ and $\Hol_{V,\omega_1}$ are symplectically (resp.\ Hamiltonian)
  isotopic then $(V,\omega_0)$ and $(V,\omega_1)$ are called isotopic (resp.\ Hamiltonian isotopic). For a fixed
  smooth fibration $V\to S^1$ the isotopy between $(V,\omega_0)$ and $(V,\omega_1)$ is equivalent to a homotopy
  of forms $\omega_0$ and $\omega_1$ through closed forms positive on fibers of the fibration.

\noindent We will prove in this section
\begin{theorem}\label{thm:filling-fibrations}
For any symplectic fibration $(V,\omega)$ over $S^1$
 there exists a compact symplectic $4$--manifold $(W,\Omega)$ with
$$\partial(W,\Omega)=(V,\omega)\,.$$
One can  additionally arrange that $H_1(W;\Z)=0$.\footnote{This was observed by Kronheimer and Mrowka,
see \cite{Kronheimer-Mrowka}.}
\end{theorem}

\noindent
The first ingredient in the proof in the following theorem of Akbulut and Ozbagci.
\begin{theorem}[See \cite{Akbulut-Ozbagci2}, Theorem 2.1]\label{thm:Akbulut}
Theorem \ref{thm:filling-fibrations} holds up to homotopy. More precise, for any symplectic fibration
$(V,\omega)$
as above there exists  a compact symplectic $4$--manifold $(W,\Omega)$ with $H_1(W)=0$  which has a structure
of a symplectic Lefschetz fibration over $D^2$ and which restricts to $S^1= \partial D^2$ as a symplectic
fibration $(V,\wt\omega)$  homotopic to $(V,\omega)$.
\end{theorem}
\noindent The proof of this theorem is based on an observation (which the authors said they learned from Ivan Smith) that
any element of the mapping class group of a closed surface can be presented as a composition of {\it
positive} Dehn twists,%
\footnote{Here is a simple argument due to Peter Kronheimer which shows this.
Take any generic genuine (i.e.\ having exceptional fibers) Lefschetz  fibration over $CP^1$ with
the fiber of prescribed genus.
Then the product of  $+1$--twists  corresponding to vanishing cycles is the identity. Therefore,
 a $-1$--twist (and hence any $-1$--twist)  is a product $+1$--twists.}
  W\,P Thurston's construction (see \cite{Thurston}) of symplectic structure on surface fibrations, and  its
adaptation by R\,E Gompf (see \cite{Gompf}) for
Lefschetz fibrations  with positive Dehn twists around exceptional fibers.
 Exploring the freedom of the construction one can  arrange  that $H_1(W;\Z)=0$.
 Indeed, for a Lefschetz fibration over $D^2$
we have
$H_1(W)\simeq H_1(F_0)/C$, where $C\subset H_1(F_0)$ is the subgroup generated by the vanishing cycles.  But the
already mentioned above fact that  the mapping class group is generated as a monoid
by positive Dehn twists allows us to make  any cycle
in $H_1(F_0)$ vanishing.

\noindent
The second ingredient in the proof of Theorem \ref{thm:filling-fibrations} is the following proposition
based on a  variation of an argument presented in \cite{Kronheimer-Mrowka}, see Lemma \ref{lm:thimble}
below.
\begin{proposition}\label{prop:flux-fixing}
Let $(W,\Omega)$ be a Lefschetz fibration over $D^2$ and $(V,\omega)$ a symplectic fibration over $S^1$ which bounds it,
$\partial(W,\Omega)=(V,\wt\omega)$.
Then for any symplectic fibration $(V,\wt\omega)$ homotopic to $(V,\omega)$ there exists a symplectic form $\wt\Omega$
on $W$ such that $\partial(W,\wt\Omega)=(V,\wt{\wt\omega})$ where $(V,\wt{\wt\omega})$ is Hamiltonian isotopic to
$(V,\omega)$.
\end{proposition}
\noindent In other words,  Proposition \ref{prop:flux-fixing} together with Theorem \ref{thm:Akbulut} imply Theorem \ref{thm:filling-fibrations}
up to {\it Hamiltonian isotopy}.
Before proceeding with the proof we recall some standard facts about the flux homomorphism.

\noindent
Let $ F_0$ be a closed oriented surface of genus $g$ with an area form $\omega$.
We denote by $\cD=\cD( F_0)$ the group of area preserving diffeomorphisms of $ F_0$ and by
$\cD_0$ its identity component. The Lie algebra of $\cD_0$ consists of symplectic vector fields, i.e.\ the vector
fields, $\omega$--dual to closed forms.  Hence,
 given an isotopy $f_t\in\mathcal{D}_0$   which connects $f_0=\Id$ with $f_1=f$ then
   for the  time-dependent vector field $v_t$ which generates $f_t$, i.e.\
 \begin{equation}\label{eq:vf}
 v_t(f_t(x))=\frac{df_t(x)}{dt},\;t\in[0,1],\;x\in  F_0\,,
 \end{equation}
  the form
 $\alpha_t=v_t\hook\omega$ is closed for all $t\in[0,1]$.
 Diffeomorphisms generated by time-dependent vector fields \eqref{eq:vf}
 dual to {\it exact} 1--forms
  form a subgroup $\cD_H$ of  Hamiltonian diffeomorphisms. This subgroup is the kernel
  of a {\it flux}, or Calabi homomorphism (see \cite{Calabi}) which is defined as follows.
  Given $v_t$ generating $f$ as in \eqref{eq:vf} as its time-one map,
  we define
\begin{equation*}
\Flux(f)=\int\limits_0^1[v_t\hook\sigma] dt\,
\end{equation*}
where $[v_t\hook\sigma]\in H^1( F_0;\R)$ is the cohomology class of the closed form $v_t\hook\sigma$.
Though $\Flux(f)$, as defined by the above formula, is independent of the choice of the path $f_t$ up to homotopy,
it may depend on the homotopy class of this path. Note, however, that when   the genus of $ F_0$ is $>1$ then
$\mathcal{D}_0$  is contractible, and hence $\Flux(f)$ is well defined as an element of $H^1( F_0;\R)$. If
$ F_0$ is the torus then $\Flux(f)$ is defined only modulo the total area of the torus, and hence it can be viewed as an element
of $H^1( F_0;\R/\Z)$.  According to \cite{Calabi},
$$\cD_H=\Ker\Flux,$$ i.e.\ two diffeomorphisms $f,g\in \mathcal{D}$ are Hamiltonian isotopic if and only if
$\Flux(f\circ g^{-1})=0$.\footnote{This can be verified as follows.
 Let $f_t\in\mathcal{D}_0,t\in[0,1]$, be any symplectic isotopy connecting
$f_0=f$ and $f_1=g$. Denote $a_t=\Flux(f_t)$ and choose a harmonic
(for some metric) $1$--form $\alpha_t$ representing $a_t\in H^1( F_0;\R)$.
Set $\beta_t=\alpha_t-\alpha_0$. By assumption we have $\beta_1=0$. Let $\varphi_t$ be the time-one map of the
symplectic flow generated by the symplectic vector field $v_t$ $\omega$--dual to $\alpha_t$. Then
for all $t\in[0,1]$ we have $\Flux(\varphi^{-1}_t
\circ f_t)=a_0$, and hence $\varphi^{-1}_t
\circ f_t$ is a Hamiltonian isotopy between $f$ and $g$.}

Therefore, Proposition \ref{prop:flux-fixing} is equivalent to
\begin{lemma}\label{lm:thimble} Suppose $(V,\omega)$ is a symplectic fibration over $S^1$
 and $(W,\Omega)$ is a symplectic Lefschetz fibration over $D^2$
such that
$\partial(W,\Omega)=(V,\omega)$ and $H_1(W)=0$. Then  for any $a\in H^1(F_0;\R)$
(or $H^1(F_0;\R/\Z)$ if $F_0$ is the torus)
there exists a symplectic form $\wt \Omega$ on $W$ such that $$\Flux(\Hol_{V,\omega}\circ\Hol^{-1}_{V,\wt\omega})=
a,$$ where
$\wt\omega=\wt\Omega|_V$.
\end{lemma}

\proof Let us recall that given an embedded path $\delta$ from a critical value $p\in D^2$
of the Lefschetz fibration to a boundary point $q\in \partial D^2$
there exists a Lagrangian disc $\Delta_\delta$, called {\it thimble}, which projects to the path $\delta$ and whose
 boundary $\partial\Delta_\delta\subset F_q=\pi^{-1}(q)$ is  a vanishing cycle. This thimble is formed by leaves
 of the characteristic foliation of the form $\Omega|_{\pi^{-1}(\delta)}$ emanating from the corresponding
 critical point.
   Let us choose disjoint embedded paths  $\delta_1,\dots,\delta_N$ from all critical values
   of the Lefschetz fibration to points inside the arc $l=[0,1/2]\subset \R/\Z=\partial D^2$.
    Let $\Delta_i=\Delta_{\delta_i}$ and $\gamma_i=\partial\Delta_i$, $i=1,\dots, N$,
be the corresponding thimbles and  vanishing cycles. Using characteristics of $\omega$ as horizontal
lines we can trivialize the fibration
$V_{1/2}=\pi^{-1}(l)\to l$.
Note that the inclusion $H_1(F_0)\to H_1(W)$ is surjective, and the kernel
of this map is generated by the vanishing cycles (independently of
  paths along which they are transported to $F_0$ from a critical point).
By the assumption, we have $H_1(W)=0$ and hence the
projections of $\gamma_i$ to $F_0$ generate $H_1(F_0)$.  Then  the cohomology classes
$D\g_i\in H^1(F_0)$, $i=1,\dots, N$, Poincar\'e dual to
$[\g_i]\in H_1(F_0)$, generate $H^1(F_0)$. In particular, we can write
$$a=\sum\limits_1^N a_iD\gamma_i.$$ Let us recall that there exists a neighborhood $U_i$
of $\Delta_i$ symplectomorphic to a disc bundle  in $T^*(\Delta_i)$. Let $(q_1,q_2,p_1,p_2)$
be the canonical coordinates in $T^*(\Delta_i)$ such that
\begin{equation}\label{eq:Omega}
\Omega|_{\Delta_i}=dp_1\wedge dq_1+dp_2\wedge dq_2,\;
 \Delta_i\subset \{p_1=p_2=0\}
 \end{equation}
 and $  U_i=\{||p||^2=p_1^2+p_2^2<\e^2\}.$
Let $\sigma\co [0,\e]\to\R$ be a non-negative function constant near $0$ and equal to $0$ near $\e$ and such that
$$\int\limits_{||p||\leq \e}\sigma(||p||)dp_1dp_2=1\,.$$
Consider a   supported in $U_i$ closed $2$--form
\begin{equation}\label{eq:eta}
\eta_i=\sigma dp_1\wedge dp_2\,.
\end{equation}
Note that the form
$$\wt\Omega=\Omega+\sum\limits_1^N a_i\eta_i$$
is symplectic, as it follows from the explicit
expressions \eqref{eq:Omega} and \eqref{eq:eta}.
Note also that the restriction of the form $\eta_i$ to the fiber containing $\gamma_i$ vanishes, and hence
for a sufficiently small $\e>0$ the form
$\wt\omega=\wt\Omega|_V$ is  positive on the fibers of the fibration $V\to S^1$.
 Let us show that
   $\Flux(\Hol_{V,\omega}\circ \Hol^{-1}_{V,\wt\omega})=a$.
   For any oriented curve
$\gamma$ in $F_0$ we have
$$\int\limits_{\Gamma}\eta_i=D\gamma_i[\gamma]\;\hbox{and}\;\omega|_{\Gamma}=0\,,$$
where $\Gamma=\g\times l$.
 Let $\wt\Gamma$ be a cylinder formed by the characteristics of $-\wt\omega$ in $V_{1/2}$
originated at $\gamma\times 1/2$, and $\wh\Gamma$ the projection of $\wt\Gamma$ to $F_0$. Note that
$\int\limits_{\wt\Gamma}\wt\omega=0$ and that the cylinders
$\wh\Gamma$ and $\wt\Gamma$ fit together into a cylinder with the same boundary as $\Gamma$.
Let us orient
$\Gamma,\wt\Gamma$
and $\wh\Gamma$ in such a way that $$\partial \Gamma=\gamma\times1/2-\gamma\times0,\;
\partial\wt\Gamma=\wt\gamma-\gamma,\; \partial\wh\Gamma=\gamma\times1/2-\wt\gamma,$$
where $\wt\gamma=\wt\Gamma\cap F_0$.
The diffeomorphism $\Hol_{V,\omega}\circ \Hol^{-1}_{V,\wt\omega}$ coincides with the projection
$F_0\to F_0\times 1/2$ followed by the holonomy along the characteristic foliation of
$-\wt\omega|_{V_{1/2}}$
Therefore,
 \begin{align}
\Flux(\Hol_{V,\omega}\circ \Hol^{-1}_{V,\wt\omega})(\gamma)=&
\int\limits_{\wh\Gamma}\omega=\int\limits_{\wh\Gamma\cup\wt\Gamma}\wt\omega\\
=&\int\limits_\Gamma\wt\omega=
\int\limits_\Gamma\sum\limits_1^Na_i\eta_i=\sum\limits_1^Na_iD\g_i(\g)=a(\gamma)\,.\tag*{\qed}
\end{align}

To finish the proof of Theorem \ref{thm:filling-fibrations} it remains to fix
 the Hamiltonian isotopy class of the holonomy diffeomorphism. This can be done using
the following standard argument from the theory of symplectic fibrations.
\begin{lemma}
\label{rm:ext-Hamiltonian} Given any Hamiltonian diffeomorphism $h\co F_0\to F_0$, consider a symplectic fibration
 $(V=F_0\times S^1,\omega)$ with $\Hol_{V,\omega}=h$. Then  there exists a symplectic form $\Omega$ on $W=F_0\times
 D^2$ such that $\partial (W,\Omega)=(V,\omega)$.
 \end{lemma}
 \begin{proof}
  Let $H_t\co F_0\to\R$, $t\in\R$, be a  $2\pi$--periodic
 time-dependant Hamiltonian whose time one map equals $h$. Suppose that $m<H_t<M$.   We can assume that $m>0$.
 Consider an embedding $f\co F_0\times S^1\to F_0\times\R^2$ given by the formula
 $$(x,t)\mapsto (x,\,\varphi=t,\,r=\sqrt{H_t(x)}),$$
 where $x\in F_0,\,t\in S^1=\R/2\pi$, and $(r,\varphi)$ are polar coordinate on $\R^2$. Let
 $\omega_0=\omega|_{F_0}$ and $\Omega_0$
 denote the split symplectic form $\omega_0+d(r^2d\varphi)$. Then $f^*\Omega_0=\omega$. On the other hand, the embedding $f$
 extends to an embedding $\wt f\co F\times D^2\to F\times\R^2$, and hence the form $\omega$ extends to a symplectic
 form
 $\Omega=\wt f^*\Omega_0$ on $D^2\times S^1$.
  \end{proof}
\noindent This finishes off the proof of Theorem \ref{thm:filling-fibrations}.
\endproof

\noindent Before proving Theorem \ref{thm:concave-filling} let us make a general remark on
gluing of symplectic manifolds along their boundaries.
\begin{remark}\label{rm:gluing}
{\rm
Let $(W_1,\Omega_1)$ and $(W_2,\Omega_2)$ be two  symplectic manifolds, and $V_1\subset\partial W_1$ and $V_2\subset
\partial W_2$ be  components of their boundaries.  Suppose we are given an orientation reversing
 diffeomorphism $f\co V_1\to - V_2$ such that $f^*\Omega_2=\Omega_1$. Then the manifold $W=W_1\mathop{\cup}\limits_{f(V_1)=V_2} W_2$
 inherits a canonical, up to  a Hamiltonian diffeomorphism, symplectic structure $\Omega$.
 Indeed, according to the symplectic neighborhood theorem the restriction $\Omega_i|_{V_i},\,i=1,2$, determines
 $\Omega_i$ on a neighborhood of $V_i$ in $W_i$ uniquely up to a symplectomorphism fixed on $V_i$.}
\end{remark}

\begin{proof}[Proof of Theorem \ref{thm:concave-filling}]
Let $(\wt W,\wt\Omega)$ be a cobordism between $(-V,\omega)$ and
a symplectic fibration $(V',\omega')$ which is provided by Theorem \ref{thm:cobord},
and $(W,\Omega)$ be a symplectic manifold bounded by $(-V',\omega')$ which we constructed in Theorem
\ref{thm:filling-fibrations}. The required cobordism  $(W',\Omega')$ we then obtain by gluing
 $(\wt W,\Omega)$ and $(W,\Omega)$ along their common
boundary, see above Remark \ref{rm:gluing}. Moreover, note that $H_1(W';\Z)=H_1(W,\Z)$.
Hence,
 one can arrange that $H_1(W';\Z)=0$.
\end{proof}

\no\begin{proof}[Proof of Corollary \ref{cor:semi}]
According to  a theorem of Giroux (see \cite{Giroux}), any contact manifold $(V,\xi)$ admits an open book decomposition.
Hence, for any symplectic form which is positive on $\xi$ we can use
Theorem  \ref{thm:concave-filling}  to find a symplectic manifold $(W,\Omega)$ with
$\partial(W,\Omega)=(-V,\omega)$.  Attaching $(W,\Omega)$ to a non-desirable component
(or components) of the boundary of a semi-filling we will transform it to a filling.
\end{proof}

\sh{An alternative proof of Theorem \ref{thm:filling-fibrations}}

The following  lemma of Kotschick and Morita (see \cite{Kotschick-Morita})
gives an alternative proof of Theorem \ref{thm:filling-fibrations}.
\begin{lemma}\label{lm:commutator}
Let $\mathcal{D}$ be the group of symplectic (i.e.\ area and orientation preserving)
 transformations of a closed surface
$( F_0,\sigma)$ where $\sigma$ is an area form with $\int\limits_{ F_0}\sigma=1$.
 Then the commutator $[\mathcal{D},\mathcal{D}]$ contains  the identity component $\mathcal{D}_0$.
 If the genus of $ F_0$ is $>2$ then the group $\mathcal{D}$ is perfect,
  i.e.\ $\mathcal{D}=[\mathcal{D},\mathcal{D}]$.
\end{lemma}
\noindent The proof of this lemma is based on Banyaga's theorem \cite{Banyaga}
 which states that $[\mathcal{D}_0,\mathcal{D}_0]=\mathcal{D}_H$, a theorem of Harer (see \cite{Harer}) that the group
$$H_1(\Gamma_g)=\Gamma_g/[\Gamma_g,\Gamma_g],$$
where $\Gamma_g$ is the mapping class group
of the surface of genus $g$,
is trivial if $g>2$ (and it is finite for $g\leq 2$), and the following formula of Lalonde and Polterovich
from \cite{Lalonde-Polterovich}.
 For any symplectomorphism $g\in\mathcal{D}$  and
  any $f\in\mathcal{D}_0$ we have $[g,f]=gfg^{-1}f^{-1}\in\mathcal{D}_0$ and
 \begin{equation*}
 \Flux([g,f])=g^*(\Flux(f))-\Flux(f)\,.
 \end{equation*}
In particular, if the linear operator $g^*\co H^1(F_0,\R)\to H^1(F_0,\R)$ has no eigenvalues $=1$ then the
formula $$f\mapsto  \Flux([g,f])$$
defines a surjective map of $\mathcal{D}_0$ onto $H^1(F_0,\R)$ (or $H^1(F_0,\R/\Z)$ if $F_0$ is the torus).
Clearly, there are a lot of diffeomorphisms $g$ with this property, and therefore
one can represent any Hamiltonian isotopy class
from $\cD_0$ as a commutator of a fixed $g\in\cD$ and a Hamiltonian diffeomorphism.

\noindent Lemma \ref{lm:commutator} allows us to extend any  symplectic fibration over a circle
whose fiber has genus $\geq 2$
 to a symplectic fibration over a surface
with boundary. The minimal genus of this surface is equal $1+m$, where $m$ is the minimal number of commutators needed
to decompose the class of $\Hol_{V,\omega}$ in the mapping class group into a product of commutators. This gives
an alternative proof of Theorem
\ref{thm:filling-fibrations} for the case when $\genus(F_0)\geq 2$. The genus restriction  is not
 a serious obstruction for  applications.
However, it is unclear whether it is possible  to improve this construction
to accommodate the condition $H_1(W)=0$.

\section{Different flavors of symplectic fillings}\label{sec:flavors}
We conclude this article by summarizing the known
relations between all existing notions of symplectic filling
which were introduced in my earlier papers.

\no  A contact manifold $(V,\xi)$ is called
 \begin{description}
 \item[(Weak)]  {\it Weakly symplectically fillable}
if there exists a symplectic manifold $(W,\omega)$ with $\partial W=V$ and  with
$\omega|_\xi>0$;
\item[(Strong)] {\it Strongly symplectically fillable}
 if there exists a symplectic manifold $(W,\omega)$ with $\partial W=V$
such that $\omega$ is exact near the boundary and there exists its primitive  $\alpha$  such that $\xi=\{\alpha|_V=0\}$
and $d\alpha|_\xi>0$;
\item[(Stein-1)] {\it Stein  (or Weinstein) fillable} if it can be filled by Weinstein symplectic manifold,
i.e.\ an exact symplectic manifold $(W,\omega)$ such that $\omega$ admits
a primitive $\alpha$ such that the Liouville vector field
$X$ which is  $\omega$--dual to $\alpha$
(i.e.\ $X\hook\omega=\alpha$) is gradient-like for a Morse function on $W$  which is constant and attains
its maximum value on the boundary.
\end{description}

\no Stein fillability admits several equivalent reformulations.
$(V,\xi)$ is Stein fillable if and only if
\begin{description}
\item[(Stein-2)] $(V,\xi)$ can be obtained by a sequence of index $1$ contact surgeries and index $2$ surgeries
along Legendrian knots with the $(-1)$--framing with respect to the framing
given by the vector field normal to the contact structure;
\item[(Stein-3)] $(V,\xi)$ is compatible with an open book decomposition which
arises on the boundary of a Lefschetz fibration
over  a disc such that the holonomy diffeomorphisms around  singular fibers are positive Dehn twists;
 \item[(Stein-4)] $(V,\xi)$ is {\it holomorphically fillable} i.e.\ there exists a complex manifold $W$ which has $V$
 as its strictly pseudo-convex
 boundary and $\xi$ is realized as the fields complex tangencies to the boundary.
 \end{description}

The equivalence of (Stein-1) and (Stein-2) follows from \cite{Eliash-Stein} or \cite{Weinstein}. The equivalence between
(Stein-2) and (Stein-3) is established in \cite{Akbulut-Ozbagci1} and \cite{Piergallini}.
The implication $\hbox{(Stein-1)}\Rightarrow\hbox{(Stein-4)}$ is established in \cite{Eliash-Stein},
 while the opposite implication
follows from \cite{Bogomolov}.

Clearly, $$\hbox{(Stein)}\Rightarrow\hbox{(Strong)}\Rightarrow\hbox{(Weak)},$$
and all these notions imply the tightness, see \cite{Eliash-filling} and \cite{Gromov-pseudo}.
As it is shown in this paper the notion of (weak/strong)
semi-fillability introduced in \cite{Eliash-Thurston} is equivalent to (weak/strong) fillability, and hence
 from now on it
should disappear.

\no Here is a summary of what is known about the relation between
 three above notions of fillability and the notion of tightness.

\medskip \noindent {\sl Tightness does not imply  weak fillability}. Such an  example was first constructed by John Etnyre
and Ko Honda in \cite{Etnyre-Honda}. More  examples were constructed by Paolo Lisca and Andr\'as Stipsicz in
\cite{Lisca-Stipsicz}.

\medskip\no{\sl Weak fillability does not imply  strong fillability.} For instance, it was shown in \cite{Eliash-torus}
that the contact  structures $\xi_n$ on the $3$--torus induced from the standard contact structure by a $n$--sheeted
covering are all weakly symplectically fillable, but not strongly fillable if $n>1$.

\medskip\no{\sl It is not known whether strong fillability implies Stein fillability}.
There are, however,  examples
of  strong symplectic fillings which are not Stein fillings.  The first example of this kind is due to
 Dusa McDuff from \cite{McDuff} who constructed   an exact symplectic manifold $(W,\omega)$
with a disconnected  contact boundary $\partial W=V_1\bigsqcup V_2$.
This manifold cannot carry a Stein structure because it is not homotopy equivalent
to a $2$--dimensional cell complex. Let us also point out that for $\dim V>3$  the notions of Stein
and strong symplectic fillability {\it do not coincide}: using a modification of the above McDuff's argument
one can construct a strongly fillable contact manifold
which cannot be a  boundary of
a manifold  homotopy equivalent to  a half-dimensional cell complex.

We will finish this section by showing that one possible notion of
fillability which seems to be intermediate between the conditions
(Weak) and (Strong) is, in fact, equivalent to strong fillability. The
Proposition
\ref{prop:exact-strong} is equivalent to Lemma 3.1 in \cite{Eliash-boundary}. It  also
appeared in
 \cite{Ohta-Ono}.
 \begin{proposition}\label{prop:exact-strong}
 Suppose that a symplectic manifold $(W,\omega)$ weakly fills a contact manifold $(V,\xi)$. Then if
 the form $\omega$ is exact near $\partial W=V$ then it can be modified into a  symplectic form $\wt\omega$
 such that $(W,\wt\omega)$ is a strong symplectic filling of $(V,\xi)$.
 \end{proposition}
 \begin{proof}
 Let $\lambda$ be a contact form  which defines $\xi$ such that $d\lambda|_\xi=\omega|_\xi$.
 According to Lemma \ref{lm:bound-behavior}  for a sufficiently small $\e>0$ and an arbitrarily large constant $C>0$
  there exists a symplectic form $\Omega$ on
 $W$ which coincides with $\omega$ outside the  $2\e$--tubular neighborhood $U_{2\e}$ of $\partial W$, and is equal to
  $$C d(t\lambda)+\omega,\; t\in[1-\e,1],$$ inside the $\e$--tubular neighborhood $U_{\e}$ of $\partial W$.
  By assumption, $\omega$ is exact near the boundary. Hence, we can assume that  $\omega=d\alpha$ in $U_\e$.
  Let $\varphi$ be a cut-off function on $U_\e$ which is equal to $0$ near $\partial W$, and is equal to $1$ near
  the other component of the boundary of $\partial U_\e$. Then if $C$ is large the form $$\wt\Omega=
  Cd(t\lambda)+d(\varphi\alpha)$$ is symplectic, and together with $\Omega$ on $W\setminus U_\e$ defines
  a strong symplectic filling
  of $(V,\xi)$.
  \end{proof}

\begin{remark}{\rm
There are known several results concerning so-called {\it concave} symplectic fillings
(which means $\partial (W,\omega)=(-V,\xi)$).
  Paolo Lisca and  Gordana Mati\v{c} proved  in \cite{Lisca-Matic}
 that any
 Stein fillable contact manifolds embeds as a separating hypersurface of contact type
 into a closed symplectic manifold (in fact a complex projective manifold).
  Selman Akbulut and Burak Ozbagci gave in \cite{Akbulut-Ozbagci2}
  a more constructive proof of this fact. Their construction
  topologically equivalent to one considered in this paper, though they did
 not considered the problem of extension of the taming symplectic form $\omega$.
  John Etnyre and Ko Honda showed
 (see \cite{Etnyre-Honda-cobord}) that any contact manifold admits a  concave symplectic filling which
 implies that a  symplectic manifold which { strongly}
 fills a contact manifold can be realized as a domain in a closed symplectic manifold.
 A different proof of this result  is given by David Gay in \cite{Gay}.
   Theorem \ref{thm:concave-filling} proven in this paper asserts
 a similar result for  weak symplectic fillings. After learning about this article John Etnyre sent  me an argument
which shows that the weak case can be deduced from the strong one, thus giving an alternative proof of
Theorem \ref{thm:concave-filling}. His idea is that by  performing
a sequence  of Legendrian contact surgeries it is possible
 to transform a contact manifold into a homology sphere and thus, taking into account
 an argument
from Proposition \ref{prop:exact-strong}, to reduce the problem  to the case considered in their paper
\cite{Etnyre-Honda-cobord} with Ko Honda.}
\end{remark}


\begin{thebibliography}{999}
 
\let\olditem\bibitem
\def\bibitem#1]#2{\olditem{#2}}


\bibitem[AO1]{Akbulut-Ozbagci1} {\bf S Akbulut}, {\bf B Ozbagci},
 {\it Lefschetz fibrations on compact Stein surfaces}, \gtref5{2001}{10}{319}{334}

 \bibitem[AO2]{Akbulut-Ozbagci2} {\bf S Akbulut}, {\bf B Ozbagci},
{\it On the topology of compact Stein surfaces},
{Int. Math. Res. Notices}, {15} (2002) 769--780

 \bibitem[Ba]{Banyaga} {\bf A Banyaga}, Sur la structure du groupe de diff\'eomorphismes qui pr\'eservent une
 forme symplectique, {Comm. Math. Helvet.} {53} (1978) 174--227

 \bibitem[Bo]{Bogomolov} {\bf F Bogomolov}, {\bf B de Oliveira},
   {\it Stein small deformations of strictly pseudoconvex surfaces},
   {Contemp. Math.} {207}, 25--41

\bibitem[Ca]{Calabi} {\bf E Calabi}, {\it On the group of
automorphisms of a symplectic manifold}, from: ``Problems in Analysis'',
Princ. Univ. Press (1970) 1--26

\bibitem[E1]{Eliash-boundary} {\bf Y Eliashberg},
{\it On symplectic manifolds which some contact properties},
 {J. of Diff. Geom.} {33} (1991) 233--238

\bibitem[E2]{Eliash-filling} {\bf Y Eliashberg},
{\it Filling by holomorphic discs and its applications}, {London
Math. Soc. Lect. Notes} {151} (1991) 45--68

\bibitem[E3]{Eliash-torus} {\bf Y Eliashberg}, {\it Unique
 holomorphically fillable contact structure on the $3$--torus},
 {Int. Math. Res. Notices}, {2} (1996) 77--82

\bibitem[E4]{Eliash-Stein} {\bf Y Eliashberg}, {\it Topological characterization of Stein manifolds of dimension
$>4$}, {International Journal of Math.} {1} (1990) 29--46

\bibitem[ET]{Eliash-Thurston} {\bf Y Eliashberg}, {\bf W\,P Thurston},
{\it Confoliations}, {University Lectures Series}, AMS, {13} (1997)

\bibitem[EH1]{Etnyre-Honda} {\bf J Etnyre}, {\bf Ko Honda}, {\it Tight
contact structures with no symplectic filling}, {Invent. Math.} 148
(2002) 609--626

\bibitem[EH2]{Etnyre-Honda-cobord} {\bf J Etnyre}, {\bf Ko Honda},
{\it On Symplectic Cobordisms}, {Math. Annalen}, {323} (2002) 31--39

\bibitem[Gay]{Gay} {\bf D Gay}, {\it Explicit concave fillings of
contact three-manifolds}, {Proc. Cam. Phil. Soc.} {133} (2002)
431--441

\bibitem[Gi]{Giroux}{\bf E Giroux}, {\it G\'eom\'etrie de contact: de
la dimension trois vers les dimensions sup\'erieures}, 
Proc. ICM-Beijing, {2} (2002) 405--414

\bibitem[Go]{Gompf} {\bf R Gompf}, {\it A topological characterization
of symplectic manifolds}, e-print \arxiv{math.SG/0210103}

\bibitem[Gr]{Gromov-pseudo} {\bf M Gromov}, {\it Pseudo-holomorphic curves
 in symplectic manifolds}, {Invent. Math.}  {82} (1985) 307--347

\bibitem[Ha]{Harer} {\bf J Harer}, {\it The second homology of the
mapping class group of an orientable surface}, {Invent. Math.} {72}
(1983) 221--239

\bibitem[KoM]{Kotschick-Morita} {\bf D Kotschick}, {\bf S Morita},
{\it Signatures of foliated surface bundles and the symplectomorphism
groups of surfaces}, \arxiv{math.SG0305182}

\bibitem[KMr]{Kron-Mrowka-monopole} {\bf P Kronheimer}, {\bf T
Mrowka}, {\it Monopoles and contact structures}, {Invent. Math.} {130}
(1997) 209--255

\bibitem[KrM]{Kronheimer-Mrowka} {\bf P Kronheimer}, {\bf T Mrowka},
{\it Witten's conjecture and Property P}, \gtref8{2004}7{295}{310}

\bibitem[LPo]{Lalonde-Polterovich} {\bf F Lalonde}, {\bf L
Polterovich}, {\it Symplectic diffeomorphisms as isometries of Hofer's
norm}, {Topology}, {36} (1997) 711--727

\bibitem[LM]{Lisca-Matic} {\bf P Lisca}, {\bf G Mati\v{c}}, {\it Tight
contact structures and Seiberg-Witten invariants}, Invent. Math.
{129} (1997) 509 -- 525

\bibitem[LS]{Lisca-Stipsicz} {\bf P Lisca}, {\bf A Stipsicz}, {\it An
infinite family of tight, not semi-fillable contact three-manifolds},
\arxiv{math.SG/0208063}

\bibitem[LPi]{Piergallini}{\bf A Loi}, {\bf R Piergallini}, {\it Compact
Stein surfaces with boundary as branched covers of $B^4$},
{Invent. Math.} {143} (2001) 325--348

\bibitem[Mc]{McDuff} {\bf D McDuff}, {\it Symplectic manifolds with
contact type boundaries}, {Invent. Math.} {103} (1991) 651--671

\bibitem[OO]{Ohta-Ono}{\bf H Ohta}, {\bf K Ono}, {\it Simple
singularities and topology of symplectically filling $4$--manifold},
{Comment. Math. Helv.}  {74} (1999) 575--590

\bibitem[OS]{Ozsvath-Szabo} {\bf P Ozsv\'ath}, {\bf Z Szab\'o}, {\it
Holomorphic discs and genus bounds}, \gtref8{2004}8{311}{334}

\bibitem[Th]{Thurston} {\bf W\,P Thurston}, {\it Some simple examples
of symplectic manifolds}, {Proc. Amer. Math. Soc.} {55} (1976)
467--468

\bibitem[We]{Weinstein} {\bf A Weinstein}, {\it Contact surgery and
symplectic handlebodies}, {Hokkaido Math. J.} {20} (1991) 241--251

\end{thebibliography}
\end{document}